\documentclass[a4paper]{amsart}

\usepackage{graphicx,color}

\theoremstyle{plain}
\newtheorem{thm}{Theorem}

\newcounter{num}
\newcommand{\name}[3]{\addtocounter{num}{1}\hbox to 1.7cm{$#1#2{#3}$\hfill} }
\newcommand{\knum}[5]{\put(10,137){\name{#1}{#2}{#3}}}

\title[Minimal grid diagrams of  13 crossing knots with arc index 13]{Minimal grid diagrams of the prime knots with  crossing number 13 and arc index 13}

\author[H. J. Lee]{Hwa Jeong Lee}
\address{Dongguk University WISE}
\email{hjwith@dongguk.ac.kr}
\author[Y. Lee]{Yoonsang Lee}
\address{Korea Science Academy of KAIST}
\email{21-075@ksa.hs.kr}
\author[C. Lee]{Chanmin Lee}
\address{Korea Science Academy of KAIST}
\email{21-083@ksa.hs.kr}
\author[Y. Park]{Yeseo Park}
\address{Korea Science Academy of KAIST}
\email{21-041@ksa.hs.kr}
\author[H. Kim]{Hun Kim}
\address{Korea Science Academy of KAIST}
\email{hunkim@ksa.kaist.ac.kr}
\author[G. T. Jin]{Gyo Taek Jin}
\address{KAIST}
\email{trefoil@kaist.ac.kr}

\date{\today}
\begin{document}

\maketitle

\begin{abstract}
We give a list of minimal grid diagrams of the 13 crossing prime nonalternating knots which have arc index 13. 
There are 9,988 prime knots with crossing number 13. Among them 4,878 are alternating and have arc index 15. Among the other nonalternating knots, 49, 399, 1,412 and 3,250 have arc index 10, 11, 12, and 13, respectively. We used the Dowker-Thistlethwaite code of the 3,250 knots provided by the program Knotscape to generate spanning trees of the corresponding knot diagrams to obtain minimal arc presentations in the form of grid diagrams.
\end{abstract}

\section{Introduction}
A \emph{grid diagram} is a knot diagram with finitely many vertical segments and the same number of  horizontal segments such that the vertical segments cross over the horizontal segments at all crossings. As shown in Figure\,\ref{fig:13n3003grid}, we draw grid diagrams without breaking underpassing arcs.

Every knot can be drawn as a grid diagram~\cite{C1995}. The \emph{arc index} of a knot $K$ is the minimum number of vertical segments (or equivalently, horizontal segments) among all grid diagrams of $K$, denoted by $\alpha(K)$. 

Imagine a vertical axis behind a grid diagram, and then replace each horizontal segment with a broken line which touches the axis at the same horizontal level. This deformation does not change the link type of the grid diagram and is an open book presentation with the number of pages equal to the number of the vertical segments of the grid diagram.  
Such an embedding is called an \emph{arc presentation}~\cite{C1995}.

\def\mythicklines{\linethickness{1.3pt}}
\def\gridthirteenC{%
\put(10,50){\line(0,1){80}}
\put(20,10){\line(0,1){110}}
\put(30,110){\line(0,1){20}}
\put(40,100){\line(0,1){20}}
\put(50,80){\line(0,1){30}}
\put(60,20){\line(0,1){80}}
\put(70,10){\line(0,1){50}}
\put(80,40){\line(0,1){30}}
\put(90,50){\line(0,1){40}}
\put(100,30){\line(0,1){30}}
\put(110,20){\line(0,1){20}}
\put(120,30){\line(0,1){50}}
\put(130,70){\line(0,1){20}}
\put(20,10){\line(1,0){50}}
\put(60,20){\line(1,0){50}}
\put(100,30){\line(1,0){20}}
\put(80,40){\line(1,0){30}}
\put(10,50){\line(1,0){80}}
\put(70,60){\line(1,0){30}}
\put(80,70){\line(1,0){50}}
\put(50,80){\line(1,0){70}}
\put(90,90){\line(1,0){40}}
\put(40,100){\line(1,0){20}}
\put(30,110){\line(1,0){20}}
\put(20,120){\line(1,0){20}}
\put(10,130){\line(1,0){20}}
}
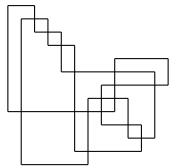
\begin{figure}[h!]
\bigskip\setlength{\unitlength}{0.5pt}
\begin{picture}(120,120)(10,10)
\gridthirteenC
\end{picture}
\caption{A minimal grid diagram of the knot $13n3003$}\label{fig:13n3003grid}
\end{figure}

\begin{thm}[\cite{Jin-Park2010}]\label{thm:less than xing}
A prime link $L$ is non-alternating if and only if
$$\alpha(L)\le c(L).$$
\end{thm}

\begin{thm}[\cite{BP2000}]\label{BP2000}
For a non-split link $L$, $\alpha(L)\le c(L)+2$.
\end{thm}

\begin{thm}[\cite{MB1998}]\label{MB1998}
Let $L$ be an alternating link, then $\alpha(L)\ge c(L)+2$.
\end{thm}

Table~\ref{table1} shows the number of prime knots of given crossing number and arc index, up to 13 crossings~\cite{N1999,Jin2020,Jin2021}.
Theorems~\ref{BP2000} and \ref{MB1998} together show that the diagonal entries of Table~\ref{table1} count alternating knots~\cite{Jin-Lee2020,Jin-Lee2022}.
There are 9,988 prime knots with crossing number 13 including 4,878  alternating ones~\cite{HTW1998}. The  alternating ones have arc index 15. Among the other non-alternating knots, 49, 399, 1,412 and 3,250 have arc index 10, 11, 12, and 13, respectively. 
Minimal grid diagrams of the 49, 399, and 1,412 knots with arc index 10, 11, and 12 appear in the articles \cite{Jin2006}, \cite{Jin-Park2011}, and \cite{Jin2020}, respectively. In this article we deal with the remaning 3,250 knots.
\begin{table}[t]
\small
\caption{Prime knots up to crossing number 13}\label{table1}
{
\newcommand{\vsp}{\vrule width0pt height 10pt depth.5pt} 
\begin{tabular}{|c|c|c|c|c|c|c|c|c|c|c|c|c|}
\hline
\vrule width0pt height 12pt depth22pt
\lower15pt\hbox{Crossings} \kern-35pt Arc index &
        \hbox to 7pt{\hfil 5\hfil}&
        \hbox to 7pt{\hfil 6\hfil}&
        \hbox to 7pt{\hfil 7\hfil}&
        \hbox to 7pt{\hfil 8\hfil}&
        \hbox to 7pt{\hfil 9\hfil}&
        \hbox to 7pt{\hfil {10}\hfil}&
        \hbox to 7pt{\hfil {11}\hfil}&
        \hbox to 7pt{\hfil {12}\hfil}&
        \hbox to 7pt{\hfil {13}\hfil}&
        \hbox to 7pt{\hfil {14}\hfil}&
        \hbox to 7pt{\hfil {15}\hfil}& Total\\
\hline
\vsp  3&1 & & & & &  &  &  &  &  &  &  1\\
\hline
\vsp  4& &1 & & & &  &  &  &  &  &  &  1\\
\hline
\vsp  5& & &2 & & &  &  &  &  &  &  &  2\\
\hline
\vsp  6& & & &3 & &  &  &  &  &  &  &  3\\
\hline
\vsp  7& & & & &7 &  &  &  &  &  &  &  7\\
\hline
\vsp  8& & &1&2& &{18} &  &  &  &  &  &21\\
\hline
\vsp  9& & & &2&6&  &{41} &  &  &  &  & 49\\
\hline
\vsp 10& & & &1&9&32&  &  123 &  &  &  &  165\\
\hline
\vsp 11& & & & &4&46&135&  & 367 &  &  & 552\\
\hline
\vsp 12& & & & &2&48&  211&   627&  &1288 &  & 2176\\
\hline
\vsp 13& & & & & &49&  399 & 1412 &  3250 & &  4878 &9988\\
\hline

\end{tabular}\\
}
\end{table}

\section{The knot-spoke method and the filtered spanning tree method}
The proof of Theorem~\ref{BP2000} introduces the knot-spoke method. It is a sequence of contractions of edges\footnote{An edge of a diagram is an arc between two neighboring crossings.} to convert a given diagram $D$ into a diagram of a single multi-crossing with  $c(D)+2$ loops which can be pushed into separate vertical half planes with a common axis to create an arc presentation. The projection of the arc presentation onto a plane perpendicular to the axis is called a \emph{wheel diagram}. It has $c(D)+2$ sticks which are called \emph{spokes}~\cite{BP2000}.

Instead of contracting edges, we attach edges one-by-one until we get a spanning tree of $D$ called a \emph{filtered spanning tree}. There are two kinds of forbidden choices for new edges; not to separate untouched crossings into two parts and except at the last step avoid choosing an edge whose one-step extension in the same direction makes a loop. The filtered spanning tree has $c(D)-1$ edges. In a vertical cylinder of a tubular neighborhood of the tree, each string is given a constant height (or depth) determined step by step according to the crossings. These heights determine the heights of the end points of the $c(D)+1$ arcs outside of the boubdary circle of the tubular neighborhood. We break the outside arc which is the one-step extension of the last edge of the tree into two part and assign the highest or the lowest height to the broken point. Contracting the tubular neighborhood to a point, we get a diagram of a single multi-crossing and $c(D)+2$ loops. From this we get a wheel diagram of $c(D)+2$ spokes~\cite{Jin-Lee2012}.

\section{Filtered spanning tree method adapted to non-alternating knots}
Every non-alternating diagram of a prime knot without any nugatory crossings has at least 4 non-alternating edges which are located on the boundaries of at least 4 regions. 
Two of the non-alternating edges can be used to reduce the arc presentation. See the proof of Theorem~\ref{thm:less than xing} in~\cite{Jin-Park2010}. 

\begin{figure}[h]
\centerline{%
\includegraphics[height=0.4\textwidth]{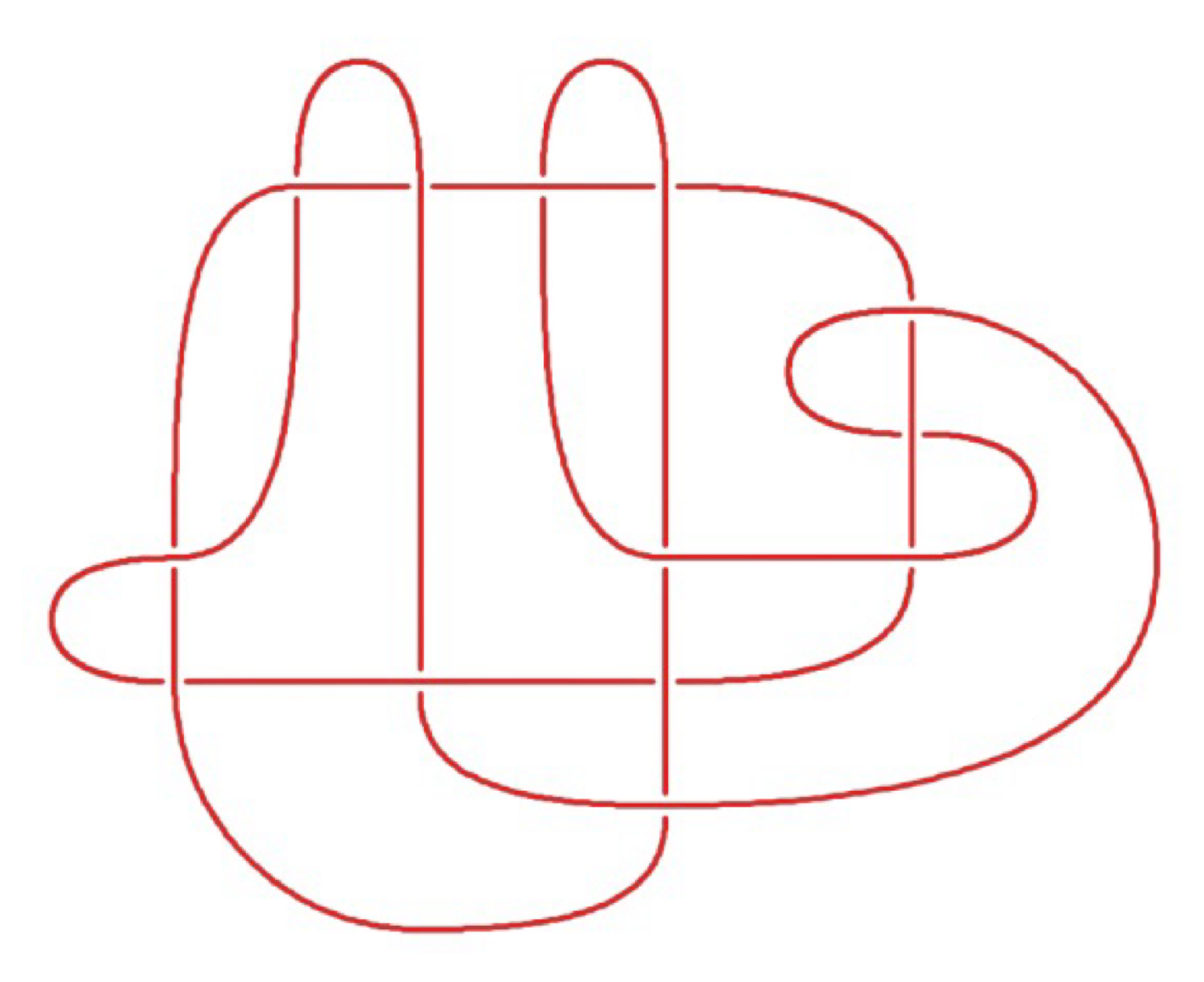}\quad
\includegraphics[height=0.4\textwidth]{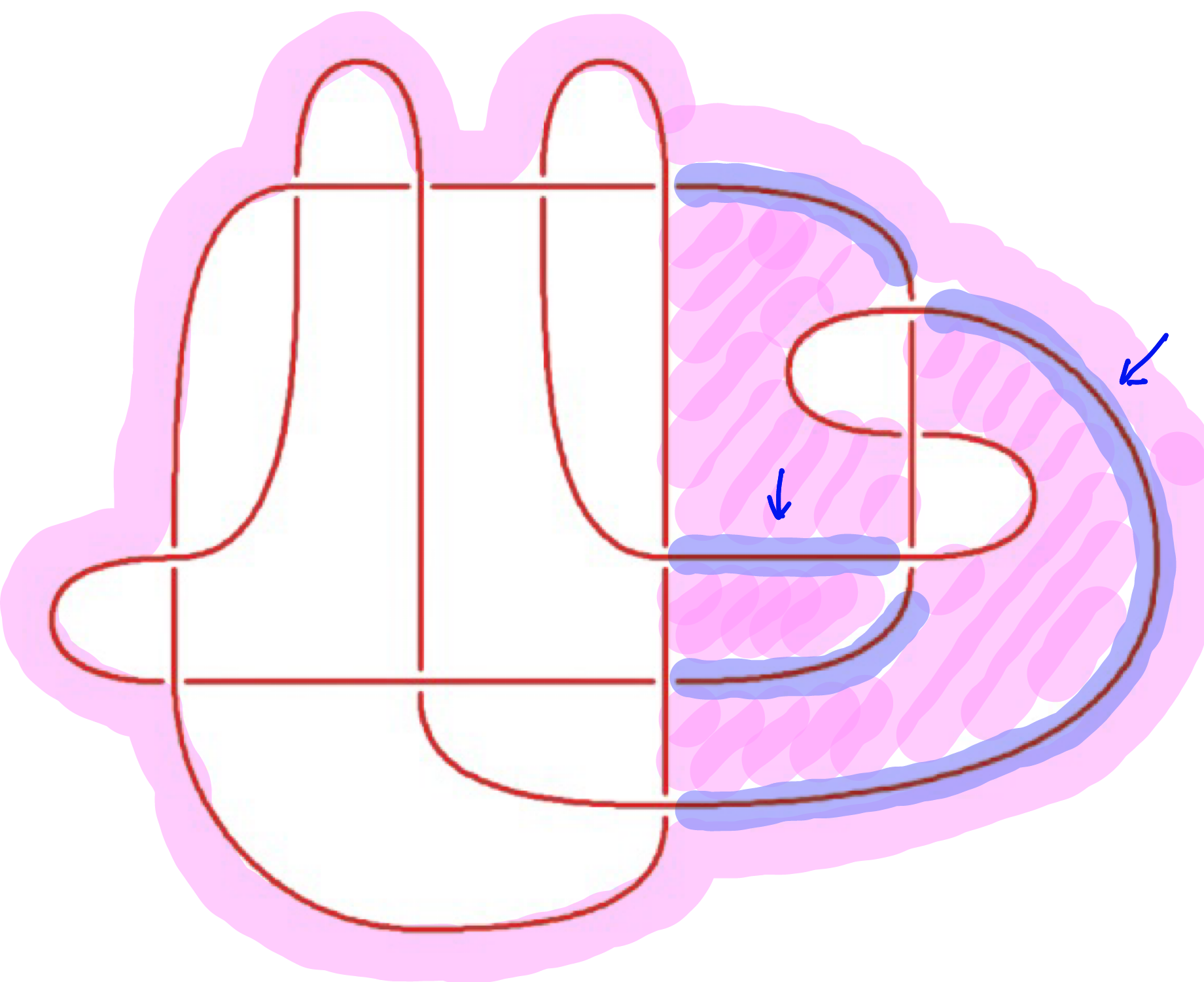}}
\caption{Diagram of $13n3003$ and its non-alternating edges \cite{SnapPy}}\label{fig:13n3003edges}
\end{figure}

\begin{figure}[h]
\centerline{%
\includegraphics[height=0.4\textwidth]{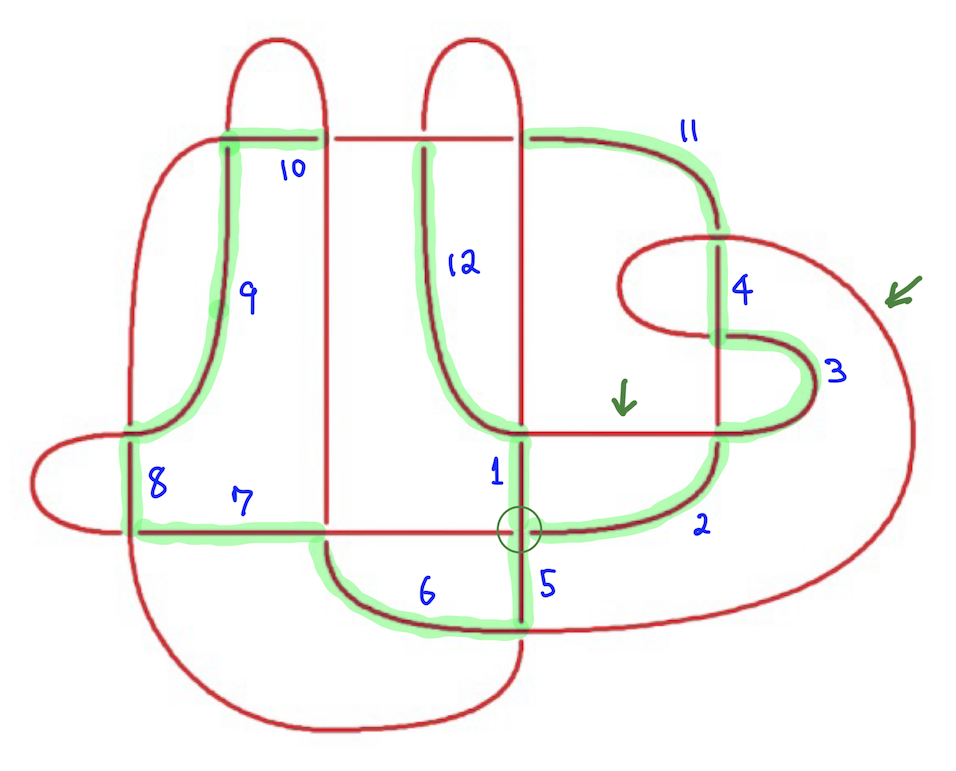}
\includegraphics[height=0.4\textwidth]{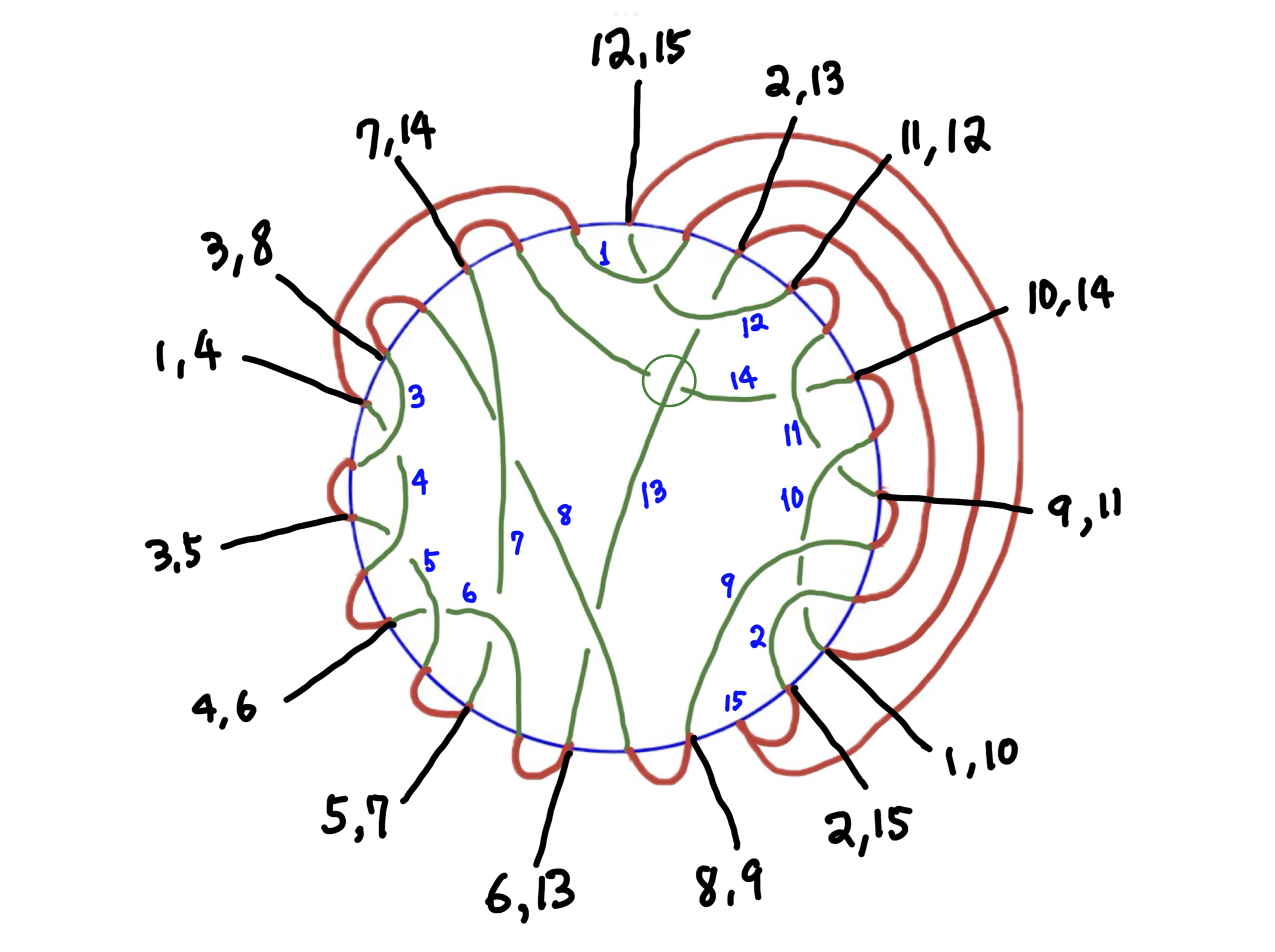}
}
\caption{A filtered spanning tree and the resulting spokes}
\label{tree-and-spokes}
\end{figure}

Figure~\ref{fig:13n3003edges} shows the knot $13n3003$. There are 4 non-alternating edges which belong to boundaries of 4 regions. 
The left of Figure~\ref{tree-and-spokes} shows a filtered spanning tree on $13n3003$ where the numbers indicate the order  of the edges of the tree. 
Contracting the edges one by one in the order as labeled is the knot-spoke method.
Notice that the edges of the tree were chosen so that each arrowed non-alternating edge is between end points of two  edges consecutively  numbered. 
The right of Figure~\ref{tree-and-spokes} is the same diagram after a planar isotopy deforming a tubular neighborhood of the tree into a circuler disk. Each arc inside the disk is at a constant depth. The arc corresponding to the edge 1 (and the edge 5) is assigned the depth 13, the edge 2 the depth 14, and so on. As one adds a new edge to the tree, the new depth is one level deeper than the maximun depth or one level shallower than the minimum depth depending on the crossing where the new edge is attached.

For each innermost arc outside the disk, we place a spoke labeled with two depths of its end points, say `$a,b$', at any one of its two end points. In the figure, they are
$$11,12;\ 10,14;\ 9,11;\ 2,15;\ 8,9;\ 6,13;\ 5,7;\ 4,6;\  3,5;\ 3,8;\ 7,14.
$$

For each arc enclosing other arcs, if one end, say $x$, has its depth between two different depths of the inner arcs and the other doesn't, then we place a spoke at $x$. If both ends are not between any two depths from inner arcs, we place a spoke at any ends. Such placement of spokes guarantees that after shrinking the cylinder of the disk the outside arcs can be pushed into the half planes lying over the spokes.
Reading the spokes of Figure~\ref{tree-and-spokes} clockwise starting from `2,13', 
we have the sequence of vertical intervals

$$\begin{aligned}
&[2,13], [11,12], [10,14], [9,11], [1,10],[2,15], [8,9], \\
&[6,13], [5,7], [4,6], [3,5], [1,4], [3,8], [7,14], [12,15].
\end{aligned}$$
In the $xy$-plane, we draw vertical segments $[2,13], [11,12], [10,14], \ldots$, on the vertical  lines $x=1,2,3,\ldots$, respectively. We add horizontal segments joining end points of the vertical segments to obtain the grid diagram on the left of 
Figure~\ref{15grid to 13grid}.

The two arrowed non-alternating edges in Figure~\ref{fig:13n3003edges}
correspond to the spokes `$11,12$' and `$8,9$'. As Figure~\ref{15grid to 13grid} illustrates, we can remove their corresponding vertical segments,  and thus obtain a grid diagram of size 13.

\def\gridthirteenA{%
\put(10,30){\line(0,1){90}}
\put(20,100){\line(0,1){30}}
\put(30,90){\line(0,1){20}}
\put(40,20){\line(0,1){80}}
\put(50,30){\line(0,1){110}}
\put(60,70){\line(0,1){50}}
\put(70,60){\line(0,1){20}}
\put(80,50){\line(0,1){20}}
\put(90,40){\line(0,1){20}}
\put(100,20){\line(0,1){30}}
\put(110,40){\line(0,1){50}}
\put(120,80){\line(0,1){50}}
\put(130,110){\line(0,1){30}}
\put(40,20){\line(1,0){60}}
\put(10,30){\line(1,0){40}}
\put(90,40){\line(1,0){20}}
\put(80,50){\line(1,0){20}}
\put(70,60){\line(1,0){20}}
\put(60,70){\line(1,0){20}}
\put(70,80){\line(1,0){50}}
\put(30,90){\line(1,0){80}}
\put(20,100){\line(1,0){20}}
\put(30,110){\line(1,0){100}}
\put(10,120){\line(1,0){50}}
\put(20,130){\line(1,0){100}}
\put(50,140){\line(1,0){80}}
}

\begin{figure}[h]
\setlength\unitlength{0.42pt}
\begin{picture}(140,140)(10,10)
\put(10,20){\line(0,1){110}}
\mythicklines\color{red}\put(20,110){\line(0,1){10}}\thinlines\color{black}
\put(30,100){\line(0,1){40}}
\put(40,90){\line(0,1){20}}
\put(50,10){\line(0,1){90}}
\put(60,20){\line(0,1){130}}
\mythicklines\color{red}\put(70,80){\line(0,1){10}}\thinlines\color{black}
\put(80,60){\line(0,1){70}}
\put(90,50){\line(0,1){20}}
\put(100,40){\line(0,1){20}}
\put(110,30){\line(0,1){20}}
\put(120,10){\line(0,1){30}}
\put(130,30){\line(0,1){50}}
\put(140,70){\line(0,1){70}}
\put(150,120){\line(0,1){30}}
\put(50,10){\line(1,0){70}}
\put(10,20){\line(1,0){50}}
\put(110,30){\line(1,0){20}}
\put(100,40){\line(1,0){20}}
\put(90,50){\line(1,0){20}}
\put(80,60){\line(1,0){20}}
\put(90,70){\line(1,0){50}}
\mythicklines\color{red}\put(70,80){\line(1,0){60}}
\put(40,90){\line(1,0){30}}\thinlines\color{black}
\put(30,100){\line(1,0){20}}
\mythicklines\color{red}\put(20,110){\line(1,0){20}}
\put(20,120){\line(1,0){130}}\thinlines\color{black}
\put(10,130){\line(1,0){70}}
\put(30,140){\line(1,0){110}}
\put(60,150){\line(1,0){90}}
\end{picture}
\begin{picture}(40,140)(10,10)
\put(20,80){$\Rightarrow$}
\end{picture}
\begin{picture}(120,140)(10,10)
\gridthirteenA
\mythicklines\color{red}
\put(30,90){\line(1,0){80}}
\put(30,110){\line(1,0){100}}
\end{picture}
\caption{Destabilization at two places}\label{15grid to 13grid}
\end{figure}
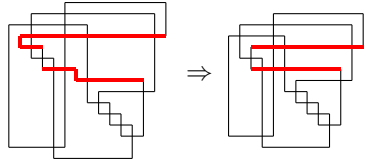

\section{Minimal Grid Diagrams of Arc Index 13}

Grid diagrams have the property that outermost edges on one side can be shifted to the other side without changing the knot type, as shown in Figure~\ref{edge moves}.
After reducing the size from 15 to 13, we applied such edge moves to make the diagrams look less random.
\def\gridthirteenB{%
\put(10,60){\line(0,1){40}}
\put(20,40){\line(0,1){30}}
\put(30,30){\line(0,1){20}}
\put(40,40){\line(0,1){50}}
\put(50,80){\line(0,1){20}}
\put(60,60){\line(0,1){80}}
\put(70,20){\line(0,1){110}}
\put(80,120){\line(0,1){20}}
\put(90,110){\line(0,1){20}}
\put(100,90){\line(0,1){30}}
\put(110,30){\line(0,1){80}}
\put(120,20){\line(0,1){50}}
\put(130,50){\line(0,1){30}}
\put(70,20){\line(1,0){50}}
\put(30,30){\line(1,0){80}}
\put(20,40){\line(1,0){20}}
\put(30,50){\line(1,0){100}}
\put(10,60){\line(1,0){50}}
\put(20,70){\line(1,0){100}}
\put(50,80){\line(1,0){80}}
\put(40,90){\line(1,0){60}}
\put(10,100){\line(1,0){40}}
\put(90,110){\line(1,0){20}}
\put(80,120){\line(1,0){20}}
\put(70,130){\line(1,0){20}}
\put(60,140){\line(1,0){20}}
}

\begin{figure}[h]
{\setlength\unitlength{0.42pt}
\begin{picture}(120,120)(10,10)
\gridthirteenA
\mythicklines\color{red}
\put(40,20){\line(1,0){60}}
\put(10,30){\line(1,0){40}}
\put(90,40){\line(1,0){20}}
\put(80,50){\line(1,0){20}}
\put(70,60){\line(1,0){20}}
\put(60,70){\line(1,0){20}}
\end{picture}
\begin{picture}(40,120)(10,10)
\put(20,70){$\Rightarrow$}
\end{picture}
\begin{picture}(120,120)(10,10)
\gridthirteenB
\mythicklines\color{red}
\put(40,90){\line(1,0){60}}
\put(10,100){\line(1,0){40}}
\put(90,110){\line(1,0){20}}
\put(80,120){\line(1,0){20}}
\put(70,130){\line(1,0){20}}
\put(60,140){\line(1,0){20}}
\end{picture}
\qquad\qquad
\begin{picture}(120,120)(10,10)
\gridthirteenB
\mythicklines\color{red}
\put(10,60){\line(0,1){40}}
\put(20,40){\line(0,1){30}}
\put(30,30){\line(0,1){20}}
\put(40,40){\line(0,1){50}}
\put(50,80){\line(0,1){20}}
\end{picture}
\begin{picture}(40,120)(10,10)
\put(20,70){$\Rightarrow$}
\end{picture}
\begin{picture}(120,120)(10,0)
\gridthirteenC
\mythicklines\color{red}
\put(90,50){\line(0,1){40}}
\put(100,30){\line(0,1){30}}
\put(110,20){\line(0,1){20}}
\put(120,30){\line(0,1){50}}
\put(130,70){\line(0,1){20}}
\end{picture}
}
\caption{Bottom 6 edges moved to the top and then left 5 edges  moved to the right}
\label{edge moves}
\end{figure}
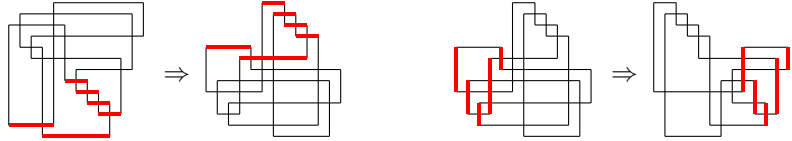

{\setlength{\unitlength}{0.4pt}
 \noindent \input k3250new2.tex 
}

\section*{Acknowledgments}
The first author was supported by  the National Research Foundation of Korea(NRF) grant funded by the Korea government(MSIT) (No. 2019R1A2C1005506).  The others' work was funded by the Ministry of Science and ICT. 

We used the program Knotscape at the beginning to get Dowker-Thistlethwaite codes of the 3,250 knots to program  for the filtered spanning trees and at the end to confirm the grid diagrams represent the original ones~\cite{knotscape}. 
We also used the program KnotPlot to obtain the Dowker-Thistlethwaite codes from the final grid diagrams~\cite{knotplot}.
The knot diagram of $13n3003$ in Figure~\ref{fig:13n3003edges} was obtained from SnapPy~\cite{SnapPy}.

\end{document}